\documentclass[12pt]{amsart}

\usepackage{amsfonts, amssymb, curves}
\newcommand{\R}{\mathbb{R}}
\newcommand{\N}{\mathbb{N}}

\newcommand{\Z}{\mathbb{Z}}

\newcommand{\T}{\mathcal{T}}

\newcommand{\A}{\mathcal{A}}

\newcommand{\cL}{\mathcal{L}}

\newcommand{\cP}{\mathcal{P}}

\newcommand{\vp}{\varphi}

\newcommand{\larr}{\left( \begin{array}{c}}
\newcommand{\rarr}{\end{array} \right) }

\newcommand{\lsqarr}{\left[ \begin{array}{c}} 
\newcommand{\rsqarr}{\end{array} \right]}

\newcommand{\arrow}{\rightarrow}

\newcommand{\inv}{\varprojlim}
\newcommand{\dir}{\varinjlim}

\begin{document}

\newtheorem{theorem}{Theorem}

\newtheorem{corollary}[theorem]{Corollary}
\newtheorem{lemma}[theorem]{Lemma} \newtheorem{prop}[theorem]{Proposition}
\newtheorem{example}[theorem]{Example}

\title{Cohomology in one-dimensional  substitution tiling
spaces}
\author{Marcy Barge and  Beverly Diamond}

 \maketitle
 \begin{abstract} Anderson and Putnam
showed that the cohomology of a substitution tiling space may be computed by
 collaring tiles to obtain a substitution
which ``forces its border.''  One can then represent the tiling space as an
inverse limit of an inflation and substitution map on a cellular complex formed
from the collared tiles; the  cohomology of the tiling space is 
computed as  the direct limit of the homomorphism induced by inflation and
substitution on the cohomology of the complex.  For one-dimensional substitution
tiling spaces, we describe a modification of the Anderson-Putnam complex on
collared tiles that allows for easier computation  and provides a means of 
identifying  certain special features of the tiling space
with particular elements of the cohomology.
\end{abstract} 


\section{Introduction}

The investigation of the topology of substitution tiling spaces, in particular,
the study of the cohomology of the space, is often aided by representing the
tiling space as the inverse limit of an inflation and substitution map on the
cellular Anderson-Putnam complex (\cite{ap}).  If the substitution {\em forces its
border}, the complex can be formed using the given tiles as cells, and adding
identifications indicating allowed transitions; the cohomology of the tiling
space is  computed as  the direct limit of the homeomorphism induced by
inflation and substitution on the cohomology of the complex.   If the
substitution does not force its border, then the inverse limit constructed in
this way may be a nontrivial quotient of the tiling space rather than the tiling
space itself;  certain arc components may be glued together, and others pinched
at points, in some cases changing the cohomology.  In this case, the process of 
{\em collaring } provides a larger set of tiles and an induced substitution which
does force its border.  The inverse limit space of the complex formed from the
larger set of tiles and associated map is homeomorphic to the tiling space of the
original substitution.

In many situations, full collaring leads to a very large number of
tiles, and the process of determining the complex can be tedious and
time-consuming.  We indicate below, for one-dimensional substitution tilings
(that is, the tiles are subintervals of $\R$),  a process which achieves the same
ultimate goal, that of representing the tiling space as an inverse limit on an
associated, simple complex, from which the cohomology can be computed in a fairly
straightforward way.

One can always {\em rewrite} a one-dimensional  substitution to obtain
a {\em proper} substitution (see \cite{dur} and \cite{bd1}) that then
automatically forces its border.  The topology of the tiling space is unchanged
by rewriting (\cite{bd1}), and the resulting
proper substitution typically has a smaller alphabet than that of the
substitution produced by collaring, hence cohomology may be more efficiently
computed by rewriting than by collaring.   This process still has two drawbacks in
our view.  First, there is no analog of rewriting in higher dimensions, while the
construction we describe below generalizes easily to higher dimensions.   Second,
the approach we present yields a specific, concrete, topological interpretation
of a piece of the cohomology.

\section{Notation and Terminology}
We introduce some notation and terminology.

Let $\A = \{1, 2, \ldots , card(\A)\}$  be a finite alphabet;   $\A^*$ will
denote the collection of finite nonempty words with letters in $\A$.  A {\em
substitution} on $\A$ is a map $\varphi : \A \arrow \A^*$; $\vp$ extends
naturally to $\vp: \A^* \arrow \A^*$.  The {\em    transition matrix} $A_\vp = A$
for $\vp$ is  $ A := (a_{ij})_{i\in
\A,  j \in \A}$ in which
$a_{ij}
$ is the number of occurrences of $i$ in the word $\vp(j)$.  

The substitution
 $\varphi$ is {\em
primitive} if
$\varphi^n(i)$ contains
$j$ for all
$i, j \in \A$ and sufficiently large $n$. Equivalently,  $\varphi$ is
primitive if and only if the matrix $A$ is aperiodic, in which case
$A$ has a simple  eigenvalue $\lambda_\varphi$ larger in modulus than its
remaining eigenvalues called the Perron-Frobenius  eigenvalue of $A$ (and
$\varphi$).

 A word $w$ is {\em allowed for $\varphi$} if
and only if for each finite subword (i.e., factor)
$w'$ of
$w$, there are $i \in \A$ and $n \in \N$ such that $w'$ is a subword
 of
$\varphi^n(i)$; the {\em language of $\vp$}, $\cL _\vp = \cL$, is the set of
finite allowed words for $\vp$. Let
$W_\varphi$ denote the set of  allowed bi-infinite words for
$\varphi$. We identify the $0^{th}$
coordinate in a bi-infinite word
$w$ by either an indexing, as in $w = \ldots w_{-1} w_0w_1 \ldots$, or by use of a
decimal point (or both). Let $\sigma: W_\vp \arrow W_\vp$ denote the {\em shift
map}:  $$ \sigma( \ldots w_{-1} . w_0 w_1 \ldots) := \ldots w_{-1}  w_0
.  w_1 \ldots .$$ The
substitution 
$\varphi : \A \arrow \A^*$ extends to $\varphi: W_\varphi \arrow W_\varphi$ 
 where
$$\varphi(\ldots w_{-1}w_0w_1 \ldots) := \ldots
\varphi(w_{-1})\ {\bf .} \ \varphi(w_0)\varphi(w_1)
\ldots.$$
The word 
$w$ is {\em  periodic} for $\varphi$ under inflation and
substitution, or {\em $\varphi$-periodic},  if for some
$m
\in
\N$, 
$$\varphi^m(w) =
\ldots
\varphi^m(w_{-1}) {\bf .} \varphi^m(w_0)\varphi^m(w_1)
\ldots = \ldots w_{-1} {\bf .} w_0w_1 \ldots.$$ 

 Each primitive
substitution $\varphi$ has at least one allowed $\varphi$-periodic bi-infinite word
which is necessarily uniformly recurrent under the shift.   A substitution $\varphi$ with precisely one periodic, hence
fixed, bi-infinite word  is called {\em proper}; $\vp$ is proper if and
only if there are
$b, e
\in
\A$ such that for all sufficiently large  $k$ 
 and all
$i \in \A$,
$\varphi^k(i) = b \ldots e$.

A primitive
substitution
$\varphi$ is {\em aperiodic} if at least one (equivalently, each)
$\varphi$-periodic  bi-infinite word is not periodic under the natural shift map,
in which case $(W_\varphi, \sigma)$ is an infinite minimal
dynamical system.  If $\vp$ is aperiodic, then the map $\varphi: W_\varphi
\arrow W_\varphi$ is one-to-one (\cite{mosse}). If
$\varphi$ is periodic (that is, primitive and not aperiodic), then
$W_\varphi$ is finite.

 Given a primitive substitution  $\varphi : \A \arrow \A^*$ with $card(\A) =
d \geq 2$, let
$\omega_L := (\omega_1, \ldots, \omega_d)$  be a positive left
eigenvector  for the Perron-Frobenius eigenvalue, 
$ \lambda$, of $A$.
 The intervals $P_i=[0,\omega_i]$,
$i=1,\ldots,d$, are called {\em prototiles for $\vp$} (consider
$P_i$ to be distinct from $P_j$ for $i\ne j$ even if $\omega_i=\omega_j$).   A
{\em tiling} $T$ of $\mathbb{R}$ by the prototiles for
$\vp$ is a collection $T=
\{T_i\}_{i=-\infty}^{\infty}$ of tiles $T_i$ for which 
$\bigcup_{i=-\infty}^{\infty}T_i=\mathbb{R}$,  each $T_i$ is a translate of
some $P_j$ (in which case we say $T_i$ is {\em of type $j$}), and $T_i\cap
T_{i+1}$ is a singleton for each
$i$.  Generally we assume that the indexing is such that $0 \in T_0 \setminus
T_{1}$.  

If $\varphi(i) = i_1 i_2 \ldots i_{k(i)}$, then
$\lambda \omega_i = \sum_{j=1}^{k(i)}\omega_{i_j}$. Thus
$|\lambda P_i| = \sum_{j=1}^{k(i)}|P_{i_j}|$, and  
$\lambda P_i$ is tiled by $\{T_j\}_{j=1}^{{k(i)}}$,
where $T_j= P_{i_j} + \sum_{k=1}^{j-1}\omega_{i_k}$. This process is called
{\em inflation and substitution} and extends to a map $\Phi$ taking a tiling
$T=\{T_i\}_{i=-\infty}^{\infty}$ of $\mathbb{R}$ by prototiles to a new tiling,
$\Phi(T)$, of $\mathbb{R}$ by  prototiles defined by inflating,
substituting,  and suitably translating each $T_i$.  More precisely, for $w
= w_1 \ldots w_n \in
\A^*$, define
$$\cP_w + t = \{P_{w_1} + t, P_{w_2} + t + |P_{w_1}|, \ldots, P_{w_n} + t +
\Sigma_{i < n}|P_{w_i}|\}.$$  Then $\Phi(P_i + t) = \cP_{\varphi(i)} +
\lambda t$ and $\Phi(\{P_{k_i} + t_i\}_{i \in \Z}) = \cup_{i
\in \Z}(\cP_{\varphi(k_i)} +
\lambda t_i)$.  

There is a natural topology on the collection $\Sigma_{\varphi}$ of all tilings
of $\mathbb{R}$ by prototiles ($\{T_i\}_{i=-\infty}^{\infty}$ and
$\{T_i^{\prime}\}_{i=-\infty}^{\infty}$ are ``close" if there is an $\epsilon$ near
$0$ so that 
$\{T_i\}_{i=-\infty}^{\infty}$ and $\{T_i^{\prime} + \epsilon\}_{i=-\infty}^{\infty}$
are identical in a large neighborhood of $0$ (see
\cite{ap} for details)).  The space
$\Sigma_{\varphi}$ is compact and  metrizable  with this topology and
$\Phi:\Sigma_{\varphi}\to\Sigma_{\varphi}$ is continuous.  
Given $T = \{T_i\}_{i = - \infty}^\infty \in \Sigma_\vp$, let
$\underline{w}(T) = \ldots w_{-1}w_0w_1 \ldots$ denote the bi-infinite word
with $w_i = j$ if and only if $T_i$ is of type $j$.   
 The  {\em
tiling space associated with}
$\varphi$, $\T_{\varphi}$, is defined as $$\T_\vp = \{ T :
\underline{w}(T) \textrm{ is allowed for } \vp\}.$$ 

 There is a natural
flow (translation) on 
$\Sigma_{\varphi}$ defined by
$(\{T_i\}_{i=-\infty}^{\infty}, t)\mapsto\{T_i - t\}_{i=-\infty}^{\infty}$. If
$\vp$ is primitive and aperiodic,
$\Phi:\T_{\varphi}\to\T_{\varphi}$ is a homeomorphism (this relies on the notion
of {\em recognizability} or invertibilty for such substitutions--see \cite{mosse}
and
\cite{sol}). Each
$T \in \T_{\varphi}$  is uniformly recurrent under the flow and has dense orbit
(i.e.,   the flow is minimal on $\T_{\varphi}$).  It follows that 
$\T_{\varphi}$ is a continuum.  

Recall that a {\em composant} of a point $x$ in a topological space $X$ is the
union of the proper compact  connected subsets of $X$ containing $x$. If $\vp$ is
a primitive substitution,  
composants and arc components in $\T_\vp$ are identical; in this case we use
the terms interchangeably.  For any  substitution
$\varphi$, the arc components of the tiling space
${\mathcal T}_\varphi$ coincide with the orbits of the natural flow (translation) on
${\mathcal T}_\varphi$.  

Tilings $T, T' \in \T_\vp$ are {\em forward asymptotic} if $\lim_{t \arrow
\infty} dist(T - t, T' - t) = 0$. Equivalently, $T =\{T_i\}_{i=-\infty}^{\infty},
T' = \{T'_i\}_{i=-\infty}^{\infty}  $ are {\em forward asymptotic} if there are
$N, M
\in
\Z$ so that $T_{N+k} = T'_{M+k} $ for all $k \geq 0$. 
Composants are {\em forward asymptotic} if they contain forward asymptotic
tilings.  Backward asymptotic tilings and composants are defined similarly.

If $f: X \to X$ is a map of a compact connected metric space $X$, then the {\em
inverse limit space} with single bonding map $f$
 is the space $$\inv f  = \{(x_0,x_1,\ldots) :
f(x_i) = x_{i-1}\;{\rm for}\; i=1,2,\ldots\}$$ with metric
$$\underline{d}(\underline{x}, \underline{y}) = \Sigma_{i \geq 0}\frac{d(x_i,
y_i)}{2^i};$$ 
$\hat{f}:\inv f\to \inv f$ will denote the natural (shift)
homeomorphism $$\hat{f}(x_0,x_1,\ldots)=(f(x_0),x_0,x_1,\ldots).$$

\section{The modified complex and its cohomology}

Given a primitive, aperiodic substitution $\vp$ on $d$ letters, the complex $K =
K_\vp $ will consist of a collection of edges representing the letters of
$\A$, another collection of edges representing allowed transitions between
letters, and certain identifications. 

  As above, $\omega_{L}= (\omega_1, \ldots, \omega_d)$ is a
left eigenvector for the Perron-Frobenius eigenvalue $\lambda $ for $\vp$.  Define
$0 < \epsilon = \min\{\frac{\omega_a}{2 \lambda}\}_{a \in \A}$.  For
$a
\in
\A$, let $$e_a := [\epsilon, \omega_a - \epsilon] \times \{a\},$$ and for $ab \in
\cL$, let $$e_{ab} := [-\epsilon,  \epsilon] \times \{ab\}.$$  Define $$K :=
[(\cup_{a
\in \A} e_a) \cup (\cup_{ab \in \cL}  e_{ab})]/ \sim ,$$ where $$(\omega_a -
\epsilon, a) \sim (-\epsilon, ab)$$ and $$(\epsilon,  b) \sim (\epsilon, ab)$$
for all $a, b \in \A$, and let $$S := \cup_{ab \in \cL}  e_{ab}/ \sim.$$
Loosely, the complex $K$ is  a wedge of $d$ circles with the branch
point `blown up' to $S$; $K/S$ is homeomorphic to the wedge of
$d$ circles.

\begin{example}\label{fib} The Fibonacci substitution: \end{example} If $\vp$ is
the Fibonacci substitution  ($\vp(1) = 12, \vp(2) = 1$), then all transitions
except 22 are allowed, and $K_\vp$ is (topologically and not to scale) shown
below.  

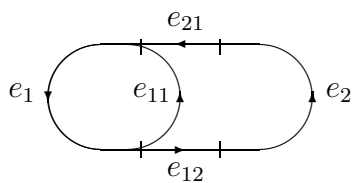
\begin{picture}(25, 80)(-100, -50)

\put(25, -10){\oval(100, 40)}
\put(-0, -10){\oval(50, 40)}
\put(7, -10){$e_{11}$}
\put(20, 17){$e_{21}$}
\put(80, -10){$e_2$}
\put(-40, -10){$e_1$}
\put(20, -40){$e_{12}$}

\put(-25, -10){\vector(0, -1){2}}
\put(75, -10){\vector(0, 1){2}}
\put(25, 10){\vector(-1, 0){2}}
\put(25, -30){\vector(1, 0){2}}
\put(25, -10){\vector(0, 1){2}}

\put(40, -33){\line(0, 1){6}}
\put(10, -33){\line(0, 1){6}}
\put(40, 6){\line(0, 1){6}}
\put(10, 6){\line(0, 1){6}}

\end{picture}

{\parindent=0pt(In }this case, collaring leads to a topologically identical
complex.) \qed

\

Given a tiling $T = \{T_n\}_{n \in \Z} \in \T_\vp$, let $t \in \R$ and  $ a, b, c
\in
\A$ be such that
$0
\in T_0 = [0, \omega_a] - t$, $T_{-1} = [0, \omega_b] - t - \omega_b$, and
$T_{1} = [0, \omega_c] - t + \omega_a$.  Define a map $p : \T_\vp \arrow K$ as
follows:  
\[ p(T) =
 \left\{
 \begin{array}{ll} [(t, a)] & \textrm{ if } 
\epsilon \leq t \leq \omega_a -
\epsilon,
\\ \left[(t, ba)\right] & \textrm{ if } 
0 \leq t \leq \epsilon,
\\
\left[(t -
\omega_a, ac)\right] & \textrm{ if }
\omega_a - \epsilon \leq t \leq
\omega_a.
\end{array} \right. \] 

{\parindent=0pt The} map $p$ is a continuous surjection, and there is a unique map
$f: K
\arrow K$ such that $f \circ p = p \circ \Phi$, where $\Phi: \T_\vp \arrow
\T_\vp$ is the inflation and substitution homeomorphism.  It is easy to verify
that if
$\vp(a) = a_1 a_2 \ldots a_k$, then $f(e_a) \subset e_{a_1} \cup  e_{a_1a_2} \cup
e_{a_2} \cup
\ldots \cup e_{a_{k-1}} \cup  e_{a_{k-1}a_k} \cup e_{a_k}$, and  if $\vp(a.b)
=
\ldots c.d \ldots$, then
$f(e_{ab})
\subseteq e_c \cup e_{cd} \cup e_d$.

\begin{lemma} $\inv f \simeq \T_\vp.$
\end{lemma}

Proof: The commuting diagram 

\begin{picture}(25, 70)(-100, -30)

\put(27, 0){$p$}
\put(90, 0){$p$}
\put(50, 27){$\Phi$}

\put(50, -18){$f$}
\put(15, -25){$K$}
\put(80, -25){$K$}

\put(15, 20){$\T_\vp$}
\put(80, 20){$\T_\vp$}
\put(22,12){\vector(0,-1){22}}
\put(85,12){\vector(0,-1){22}}

\put(75,-23){\vector(-1,0){43}}
\put(75,22){\vector(-1,0){43}}

\end{picture}

\noindent induces a map $\hat{p}: \T_\vp \arrow \inv f$
given  by $\hat{p}(T) = (p(T), p(\Phi^{-1}(T)),  \ldots )$.

 Let $V := \{ (\epsilon,
a)\}_{a \in \A} \cup \{ (\omega_a - \epsilon,
a)\}_{a \in \A} $, and choose $\underline{x} = (x_0, x_1, \ldots) \in \inv
f$. Suppose that for some
$i
\in \N$,  
$x_i\in V$.  The choice of $\epsilon$
implies that
$x_{i-1}
\in e_b \setminus V$ for some $b \in \A$.  That is, if $
T \in \hat{p}^{-1}(\underline{x}) $, then the
$0^{th}$ tile of
$\Phi^{1-i}(T)$, and the position of the origin within the interior of this tile,
are determined, as are segments of
$T $ of approximate length $\lambda^i \omega_b$
about the origin.
 It follows that 
$T \in \hat{p}^{-1}(\underline{x}) $ is totally determined if $x_i \in V$ for
arbitrarily large $i$. 

 If
$x_i \notin V$ for arbitrarily large 
$i$, then for sufficiently large  $n \in \N$, either $x_n \in e_{ab}$ for some $a,
b
\in
\A$ or
$x_n \in e_b \setminus V$ for some $b \in \A$.  In the first case,  the
$(-1)^{st}$ and $0^{th}$ tiles of
$\Phi^{1-n}(T)$, and the position of the origin within the interior of the union
of these edges, are determined.  In both cases, arguments identical to those
above imply that arbitrarily large segments of  $T$ about the origin are
determined, thus $T$ is totally determined, and
$\hat{p}$ is one-to-one.  As
$p$ is onto, so is $\hat{p}$. \qed

\

There are two essential ways in which the cohomology of the
inverse limit of the natural map  on the wedge of circles, that is, $ \dir
A_\vp^t$,  may need to be modified to obtain  the cohomology of the tiling space
$\T_\vp$.

First, if $\T_\vp$ contains one or more cycles of asymptotic composants associated
with a cycle of $\vp$-periodic words of the
form
$\ldots a_1.a_2 \ldots$, $\ldots a_3.a_2 \ldots$, $ \ldots a_3.a_4 \ldots $,
$\ldots$,
$\ldots a_1.a_{2n} \ldots$,  those composants are `glued together' in the 
inverse limit on the wedge of circles to form
$n$-ods (possibly overlapping).  The modified collaring described above, in
particular, the construction of the subcomplex $S$, `pulls apart' these cycles. 
The cohomology associated with the unglued  cycles in
$\T_\vp$ is computed via
$\check{H}^1(S)$ and must be added to $\dir A^t$.

Second, there may be `extra'  cohomology appearing in $\dir A^t$ that does not
appear in
$\T_\vp$ that we identify below.  

We find it easier to compute the
cohomology of
$\T_\vp$ by working with a map
$g$ defined on
$K$ that is  homotopic to
$f$ and leaves the transition edges invariant.  We define $g$ as follows.

First, let $\delta = \min\{\frac{\omega_a}{3}\}_{a \in \A}$.  For $a \in
\A$, define
\[ h_a(x)  =
 \left\{
 \begin{array}{ll} \epsilon & \textrm{ if } 
\epsilon \leq x \leq \delta,
\\ \epsilon + \frac{\omega_a - 2 \epsilon}{\omega_a - 2 \delta}(x-\delta) &
\textrm{ if } 
\delta \leq x \leq \omega_a - \delta,
\\

\omega_a - \epsilon  & \textrm{ if }
\omega_a - \delta \leq x \leq
\omega_a - \epsilon.
\end{array} \right. \] 

{\parindent=0pt That }is, $h_a$ linearly maps the interval $[\delta, \omega_a -
\delta]$ over
$[\epsilon, \omega_a - \epsilon]$ and collapses the intervals $[\epsilon,
\delta]$ and $[ \omega_a - \delta,  \omega_a - \epsilon]$ to 
$\{\epsilon\}$ and $\{\lambda_a- \epsilon\}$, respectively.

Let $h : K \arrow K$ be defined as:
\[ h(x, w)  =
 \left\{
 \begin{array}{ll} (h_a(x), w) & \textrm{ if } 
 x \in e_a, \textrm{ for } a \in \A,
\\ (x, w) &
\textrm{ if } 
x \in e_{ab}, \textrm{ for }  ab \in \cL.
\end{array} \right. \]

 Finally, $g := h \circ f$.  If $\vp(a.b) = \ldots c.d \ldots$, then $g(e_{ab}) =
e_{cd}$, hence $g(S) \subseteq S$.  

 We show below that the `extra' piece of cohomology mentioned  above may be
computed 
 as follows.   We find the eventual range
$ER$ of
$g|_S$, and the rank of the group generated by the coboundaries of the components
of
$ER$.  Each generator for this group is  an eigenvector associated with
eigenvalue  1 for a power of $g^*$ acting on  $H^1(K, S)$, and the rank
of the group indicates the number of copies of $\Z$ to be quotiented out from
$\dir A^t$.

Whereas $\inv f \simeq \T_\vp$,  $\inv g$ is homeomorphic to a quotient of
$\T_\vp$ in which asymptotic composants for $\vp$ associated with $\vp$-periodic 
words, should there be any,  are glued in the direction of asymptoticity (each
such pair of composants is associated with a pair of $\vp$-periodic  words of the
form
$\ldots a.b
\ldots
$,
$\ldots a.c\ldots
$, thus with the edges $e_a$,  $e_{ab}$ and $e_{ac}$, with the last two periodic
under $g$).  However, in this case, unlike the situation with the inverse
limit space on the wedge of circles, gluing does not extend all the way to the
origins of the associated tilings, and the next result says that this gluing does
not change the cohomology of
$\T_\vp$.

\begin{corollary} $\dir g^* \simeq \check{H}(\T_\vp)$.
\end{corollary}

Proof:  Since $f$ and $g$ are homotopic, $\dir g^* = \dir f^*
\simeq
\check{H}(\T_\vp)$. \qed

\

Denote the eventual range $\cap_{n \geq 0} g^n(S)$ of $g|_S$ by  $ER$.  Let
$k$ be the number of connected components of $ER$ and $l$ the number of
independent 1-cycles in $ER$; we call $l$ the {\em asymptotic cycle rank of
$\vp$}. (Note that if
$\vp$ is proper, $ER$ consists of a single edge $e_{ab}$. If $\vp$ forces its
border, then each component of $ER$ is a single edge.)

Our main result is the following:

\begin{theorem} \label{main} $ \check{H}^1(\T_\vp) \simeq   \Z^l
\oplus
\dir A^t/G 
$, where $G \simeq \Z^{k-1}$.
\end{theorem}

Proof:  There is an
$m
\in
\N$ such that $g^m(C) = C$ for each component $C$ of $ER$ and $g^m(S) \subseteq
ER$. Pick a point $c$ from each component $C$ of $S$, and let $\hat{c}$ denote the
corresponding dual generator in $H^0(S)$.  For each component $C$ of $ER$, the
element
$$v_C =
\Sigma\{
\hat{c}': C' \textrm{ is a component of }S, g^m(C') \subseteq C\}$$ is fixed by
$(g^*)^m: H^0(S) \arrow H^0(S) $.  Furthermore, the eventual range of $g^*$ has
basis $\{v_C: C \textrm{ a component of } ER\}$, and the sum of these basis
elements, $\Sigma v_C$, 
generates the
range $\eta$  of the augmentation map.  Thus, choosing $k-1$ of the components,
$C_1, \ldots, C_{k-1}$, of $ER$, the cosets $\{v_{C_1} + \eta, \ldots,
v_{C_{k-1}} + \eta\}$ form a basis for the eventual range of $g^*$ on the reduced
cohomology $\tilde{H}^0(S) := H^0(S) /\eta$, and each of these basis elements is
an eigenvector  for $(g^*)^m$ with eigenvalue 1.

It is clear that $g^m$ is homotopic to the identity on $ER$,  so that $\dir g^*:
H^1(S) \arrow H^1(S) \simeq \Z^l$, where $l$ is the asymptotic cycle rank of
$\vp$.  To see how the 1-cycles  of $ER$ contribute to the cohomology of
$\T_\vp$, consider the morphism of the exact sequence of the pair
$(K, S)$ induced by $g:(K, S) \arrow (K, S)$:

{\parindent=0pt
\begin{picture}(20, 75)(-115, -35)

\put(-50, 20){$0 \arrow \tilde{H}^0(S) \stackrel{\delta}\arrow {H}^1(K, S)
\arrow {H}^1(K ) \arrow {H}^1(S ) \arrow 0$
}

\put(-50,-20){$0 \arrow \tilde{H}^0(S) \stackrel{\delta}\arrow {H}^1(K, S)
\arrow {H}^1(K ) \arrow {H}^1(S ) \arrow 0$
}

\put(-10,15){\vector(0,-1){20}}

\put(50,15){\vector(0,-1){20}}
\put(110,15){\vector(0,-1){20}}
\put(160,15){\vector(0,-1){20}}

\put(-5, 3){$g^*_0$}
\put(55, 3){$g^*_1$}
\put(115, 3){$g^*_2$}
\put(165, 3){$g^*_3$}

\end{picture}}

{\parindent=0pt Taking }direct limits, we get an exact sequence $$0 \arrow G_0
\stackrel{\vec{\delta}}\arrow G_1
\arrow G_2 \arrow G_3  \arrow 0$$ in which $G_1 \simeq \dir A^t$,  $G_2 := \dir
g^*_2
\simeq
\check{H}^1(\T_\vp)$, $G_3 := \dir g^*_3 \simeq \Z^l$, and $G_0 := \dir g^*_0
\simeq \Z^{k-1}$ as described above.  The short exact sequence $$0 \arrow
coker(\vec{\delta}) 
\arrow G_2
 \arrow G_3  \arrow 0$$ splits to give $\check{H}^1(\T_\vp) \simeq
coker(\vec{\delta})
\oplus \Z^l$. 

 Since $H^1(K)$ is free, so is $coker(\delta)$.  In particular,
the image under $\delta$ of a basis of $\tilde{H}^0(S)$, $$\{ \delta(\hat{c}_1 +
\eta),
\ldots ,
\delta(\hat{c}_{k-1} +
\eta)\} \cup \{ \delta(\hat{c} + \eta): C \textrm{ a component of } S, C
\cap ER = \emptyset\},
$$ extends to a basis for
$H^1(K, S)$ with respect to which
$(g_1^*)^m$ takes the form $  \left( \begin{array}{ccc} I &0 &E_2
\\ E_1 & 0 & E_3
\\ 0 & 0 &A_1
\end{array} \right), $
where $I$ is the $(k-1) \times (k-1)$ identity matrix, and the diagonal  $0$
matrix is  size $j \times j$, where $j := (\# $ of components in $S) - k $. 
Note that in the basis $\{\hat{e}_a: a \in \A\}$, $g^*_1$ is represented by
$A^t$, the transpose of the transition matrix for $\vp$.  Thus $  \left( \begin{array}{ccc} I &0 &E_2
\\ E_1 & 0 & E_3
\\ 0 & 0 &A_1
\end{array} \right)  $ is unimodularly equivalent
to $(A_\vp^t)^m$, and $coker(\vec{\delta}) \simeq \dir (\overline{g^*_1})^m$,
where
$\overline{g^*_1} : coker(\delta) \arrow coker(\delta)$ is given by 
$\overline{g^*_1}(x + \delta(\tilde{H}^0(S)) = g^*_1(x) +
\delta(\tilde{H}^0(S))$.  Thus
$coker(\vec{\delta})
\simeq \dir A_1$ and we have $ \check{H}^1(\T_\vp) \simeq \dir A_1 \oplus \Z^l$.
\qed

Note:  $\dir \left( \begin{array}{cc} I & 0 
\\ E_1 & 0 
\end{array} \right)  \simeq \Z^{k-1}$.  There are conditions that imply that
$\dir A^t$ splits as $\Z^{k-1} \oplus \dir A_1$.  For instance, if there is a
matrix $B$ so that $B A_1 + \left( \begin{array}{cc} I & 0 
\\ E_1 & 0 
\end{array} \right) = - \left( \begin{array}{c} E_1
\\ E_2
\end{array} \right)$, then $\dir A^t$ splits as described.  We do not know if
there is always such a splitting.

We rephrase the above discussion in a way that is more consistent with actual
computations (for instance, see Example \ref{disS}) and allows for a more
precise statement of Theorem \ref{main}.  

Let $f_i = (0 \ldots 0 1 0 \ldots )^t$ be the $i^{th}$ standard basis vector for
$\R^d$ and $C_1, \ldots, C_p$  the components  of the complex $S$.  For $i =
1, \ldots, p-1$, let $w_i := \Sigma f_j - \Sigma f_m$, where the first sum is
over those $j$ for which $e_j = [\epsilon, \lambda_j - \epsilon] \times \{j\}$
terminates in $C_i$ (i.e., $(\lambda_j - \epsilon, j) \in C_i$) and the second 
sum is over those $m$ for which $e_m $ originates in $C_i$ (i.e., $(
\epsilon, m) \in C_i$).  Then $\{ w_1, \ldots, w_{p-1}\}$ extends to a basis $\{
w_1, \ldots, w_{p-1}, w_p, \ldots, w_d\}$ of $\Z^d$.  Let $P$ be the matrix with
columns $\{w_1, \ldots, w_{p-1}, w_p, \ldots, w_d\}$.  Then $P^{-1} A^t P$ has
the form $  \left( \begin{array}{cc} E & F
\\ 0 & A_1
\end{array} \right)  $ in which $A_1$ is size $(d - p + 1) \times (d - p + 1)$.

\

{\parindent=0pt {\bf Theorem 4$'$}.  Let $\vp$} be an aperiodic primitive
substitution with asymptotic cycle rank $l$ and  $A_1$ as
above.  Then  $ \check{H}^1(\T_\vp) \simeq \dir A_1 \oplus \Z^l$.

\begin{corollary}  Let $\vp$, $K$, and $S$ be as above.  If $\tilde{H}^*(S)
\simeq 0$, then  $ \check{H}^1(\T_\vp) \simeq  
\dir A_\vp^t $.
\end{corollary}


\begin{example} The Fibonacci substitution: \end{example} If $\vp$ is the
Fibonacci substitution  ($\vp(1) = 12, \vp(2) = 1$) (see Example \ref{fib}), then
$S$ has a single connected component, and the eventual range of $g|_S$ is the
connected pair of edges $ER = e_{11}  \cup e_{21}$.  In this case,
$\tilde{H}^0(S) \simeq 0 \simeq H^1(S)$, and $\check{H}^1(\T_\vp) \simeq \dir
\left[ \left ( \begin{array}{cc}
1 & 1 
\\ 1 & 0
\end{array} \right)^t \right] \simeq \Z^2 $.

\begin{example} The Morse-Thue substitution: \end{example} Let $\vp$ denote the
Morse-Thue substitution  ($\vp(1) = 12, \vp(2) = 21$).  

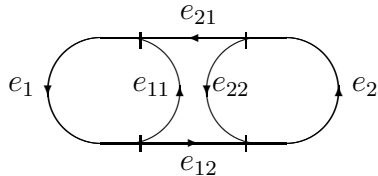
\begin{picture}(25, 75)(-100, -45)

\put(25, -10){\oval(110, 40)}
\put(-5, -10){\oval(50, 40)}
\put(55, -10){\oval(50, 40)}
\put(2, -10){$e_{11}$}
\put(32, -10){$e_{22}$}
\put(20, 17){$e_{21}$}
\put(85, -10){$e_2$}
\put(-45, -10){$e_1$}
\put(20, -40){$e_{12}$}

\put(-30, -10){\vector(0, -1){2}}
\put(80, -10){\vector(0, 1){2}}
\put(25, 10){\vector(-1, 0){2}}
\put(30, -10){\vector(0, -1){2}}
\put(25, -30){\vector(1, 0){2}}
\put(20, -10){\vector(0, 1){2}}

\put(45, -33){\line(0, 1){6}}
\put(5, -33){\line(0, 1){6}}
\put(45, 6){\line(0, 1){6}}
\put(5, 6){\line(0, 1){6}}

\end{picture}

In this case, all
transitions are allowed, and $K$ is as above.

{\parindent=0pt $S$ has a single }connected component, $ ER = S$, $\tilde{H}^0(S)
\simeq 0
$ and $H^1(S) = \Z$. Then $\check{H}^1(\T_\vp) \simeq \dir \left(
\begin{array}{cc} 1 & 1 
\\ 1 & 1
\end{array} \right)  \oplus \Z  \simeq \Z[1/2] \oplus \Z$.

\begin{example}Disconnected $S$: \label{disS}
\end{example}  Let $\vp$ be the substitution: $\vp(1) = 12341$, $\vp(2) = 12$,
$\vp(3) = 3423$, $\vp(4) = 42$.  The set $S$ has two connected
components associated with the two sets of allowed transitions: $\{ 11, 41,
21, 23, 12, 42\}$ and $\{34\}$.  The set $ER$ consists of two connected
components associated with $\{21, 11, 23\}$ and $\{34\}$.  Then $\tilde{H}^0(S)
\simeq \Z$ has basis $\{\hat{c}_1\}$, where $\hat{c_1}$ is the dual generator
of
 a vertex $c_1$ in
$e_{34}$.  It follows that 
$\{\delta \hat{c}_1 = \hat{e}_3 - \hat{e}_4\}$ extends to the basis $\{ \hat{e}_3
- \hat{e}_4, \hat{e}_1, \hat{e}_2, \hat{e}_3 \}$ of $H^1(K, S)$ with respect to
which $g_1^*$ has matrix $ \left(
\begin{array}{c ccc} 1 & 0&-1 &0 
\\ 0 & 2 & 1 & 1\\0 & 1 & 1 & 0\\ 0 & 0 & 2 &2
\end{array} \right)$.  Since there are no 1-cycles in $S$, 
$\check{H}^1(\T_\vp)  \simeq \dir  \left(
\begin{array}{ccc}  2 & 1 & 1\\ 1 & 1 & 0\\  0 & 2 &2
\end{array} \right)$.

{\Small
{\parindent=0pt
Department of Mathematics, Montana State University, Bozeman, MT 59717

barge@math.montana.edu

Department of Mathematics, College of Charleston, Charleston, SC 29424
 
diamondb@cofc.edu}}

\end{document}